\documentclass[11pt]{article}  

\usepackage{amsmath, amsthm, amsfonts, amssymb}
\usepackage[bookmarks]{hyperref}
\usepackage{indentfirst} 
\usepackage{marvosym}

\usepackage[a4paper]{geometry}

\newtheorem{lemma}{Lemma}[section]
\newtheorem{theorem}[lemma]{Theorem}

\newtheorem{corollary}[lemma]{Corollary}

\theoremstyle{definition}

{\bf}{\rm}

\newcommand{\op}{\operatorname}
\newcommand{\ce}[2]{\pmb{\op{C}}_{#1}(#2)}
\newcommand{\no}[2]{\pmb{\op{N}}_{#1}(#2)}
\newcommand{\ze}[1]{\pmb{\op{Z}}(#1)}

\newcommand{\fit}[1]{\pmb{\op{F}}(#1)}

\newcommand{\syl}[2]{\op{Syl}_{#1}\left(#2\right)}
\newcommand{\hall}[2]{\op{Hall}_{#1}\left(#2\right)}

\newcommand{\irr}[1]{{\op{Irr}}(#1)}

\begin{document}

\title{\bf Finite groups and the class-size prime graph revisited}

\author{\sc  V. Sotomayor
\thanks{Instituto Universitario de Matemática Pura y Aplicada (IUMPA-UPV), Universitat Polit\`ecnica de Val\`encia, Camino de Vera s/n, 46022 Valencia, Spain. \newline
\Letter: \texttt{vsotomayor@mat.upv.es} \newline \rule{6cm}{0.1mm}\newline
This research is supported by Proyecto PGC2018-096872-B-I00 (MCIU / AEI / FEDER, UE), and by Proyectos AICO/2020/298 and CIAICO/2021/163, Generalitat Valenciana (Spain). Additionally, the author would like to thank S. Dolfi, E. Pacifici and L. Sanus for introducing him to the research on the class-size prime graph. \newline
}}

\date{}

\maketitle

\begin{abstract}
\noindent We report on recent progress concerning the relationship that exists between the algebraic structure of a finite group and certain features of its class-size prime graph.
\end{abstract}


\section{Introduction}

Hereafter, all groups considered are finite. In recent decades, the influence on the structure of a finite group $G$ of the arithmetical properties of its set $cs(G)$ of conjugacy class sizes has been a widely investigated area. In order to better understand the arithmetical structure of that set of positive integers, the use of the prime graph built on the class sizes has been gaining an increasing interest. In general, the prime graph $\Delta(X)$ built on a set of positive integers $X$ is the (simple undirected) graph whose vertex set consists of the prime divisors of the numbers in $X$, and the edge set contains pairs $\{p,q\}$ of vertices such that the product $pq$ divides some number in $X$. In particular, when $X=cs(G)$, we briefly write $\Delta(G)$ for the class-size prime graph, $V(G)$ for the vertex set, and $E(G)$ for the edge set.

The relationship that exists between graph-theoretical properties of $\Delta(G)$ and the structure of $G$ itself has attracted the interest of many authors. In this setting, two main questions naturally arise: what can be said about the structure of $G$ if some information on $\Delta(G)$ is known, and which graphs can occur as $\Delta(G)$ for some group $G$? We usually refer to the expository paper \cite{CCsurvey} due to Camina and Camina for an overview of results on this framework. Besides, the survey \cite{Lewis} of Lewis, which is mainly devoted to summarise results regarding graphs associated with character degrees, has also a section focused on graphs built on conjugacy class sizes. Nevertheless, the aforementioned papers were publised in 2011 and 2008, respectively, and numerous interesting results have appeared since then. Thus, our aim in this note is to bring them together in one place, providing references to the literature for their proofs. For the sake of completeness, the corresponding results that were already gathered in \cite{CCsurvey} or \cite{Lewis} will be again presented here.

The structure of the paper is mostly based on the two above mentioned questions. More concretely, in Section \ref{structure} we focus on how graph-theoretical properties of $\Delta(G)$ are reflected and influenced by the algebraic structure of $G$. On the other hand, in Section \ref{possible_graphs} we pay attention to the graphs that may appear as $\Delta(G)$ for some  group $G$. In the last section, we review variations of these results, which consider the class-size prime graph built on some smaller subsets of elements of the group, as $p$\emph{-regular} elements, $p$\emph{-singular} elements, \emph{real} elements, \emph{vanishing} elements, elements lying in a normal subgroup, and elements that belong to the factors of a factorised group.

\medskip

\noindent \textbf{Notation:} for a finite group $G$, let $x^G=\{g^{-1}xg \; : \; g\in G\}$ be the conjugacy class of an element $x\in G$, and let $|x^G|$ be its size. Recall that, in virtue of the orbit-stabiliser theorem, $|x^G|=|G:\ce{G}{x}|$. As previously said, $cs(G)=\{ |x^G| \; : \; x\in G\}$. By $\operatorname{Irr}(G)$ we mean the set of all irreducible complex characters of the group $G$. The sets of all Sylow $p$-subgroups and Hall $\pi$-subgroups of $G$ are $\syl{p}{G}$ and $\hall{\pi}{G}$, respectively, where $p$ is a prime number and $\pi$ is a set of primes. The set of all prime divisors of the order of $G$ is $\pi(G)$. The semidirect product of $A$ and $B$, where $A$ is normalised by $B$, is written by $A\rtimes B$. Finally, we use CFSG to denote the classification of finite simple groups. The remaining notation and terminology is standard within Finite Group Theory.

\section{The interplay between $\Delta(G)$ and the structure of $G$}
\label{structure}

Perhaps the earliest result in this strand of research is the following one due to Sylow in 1872: if $cs(G)$ only contains powers of a fixed prime $p$, then the centre of $G$ is non-trivial. This fact, that elementarily follows from the class equation, can actually be rephrased in terms of the class-size prime graph:

\begin{lemma}\emph{(\cite{Sylow}, rephrased).}
Let $G$ be a group such that $\Delta(G)$ is an isolated vertex. Then $\ze{G}$ is non-trivial.
\end{lemma}

In fact, Baer in 1953 proved that all the class sizes of prime power order elements of a group are powers of a fixed prime $p$ if and only if $G$ has a central $p$-complement (cf. \cite[Proposition 1]{B}). Thus, we have the next corollary.

\begin{corollary}
Let $G$ be a group. Then $\Delta(G)$ is an isolated vertex $p$ if and only if $G=P\times H$ with $P\in\syl{p}{G}$ and $H\leqslant\ze{G}$.
\end{corollary}

For the ``dual'' condition on the class sizes, i.e. a fixed prime $p$ does not divide any class size of a group $G$, it is not difficult to prove that  this happens if and only if $G$ has a central Sylow $p$-subgroup. In other words, we have a characterisation of the absence of a vertex in the class-size prime graph:

\begin{lemma}
\label{vertex_out}
Let $G$ be a group, $p\in\pi(G)$ and $P\in\syl{p}{G}$. Then $p$ does not lie in $V(G)$ if and only if $P\leqslant \ze{G}$.
\end{lemma}

In this spirit, the absence of an edge in $\Delta(G)$ has also a strong impact on the group structure. This fact was studied in 1953 by Itô, but without any mention of a graph. More concretely, he proved in \cite[Proposition 5.1]{Ito} that whenever $p\neq q$ are primes that divide two distinct conjugacy class sizes of a group $G$, and $pq$ does not divide any class size of $G$, then such a group has either a normal $p$-complement or a normal $q$-complement. It also follows from Itô's proof that the Sylow $p$-subgroups (or $q$-subgroups) are abelian, although he did not explicitly stated it. Indeed, Dolfi obtained in \cite[Theorem 13]{dolfi1995} that both the Sylow $p$-subgroups and Sylow $q$-subgroups are also abelian if we assume the solubility of the group. However, Casolo and Dolfi moved beyond that, obtaining statement (ii) of the following theorem.

\begin{theorem}
\label{ito}
Let $G$ be a finite group. Suppose that $p\neq q$ are non-adjacent vertices of $\Delta(G)$. Then the following conclusions hold:
\begin{itemize}
\setlength{\itemsep}{-1mm}
\item[\emph{(a)}] \emph{(\cite[Proposition 5.1]{Ito}).} $G$ has, up to interchanging $p$ and $q$, a normal $p$-complement and abelian Sylow $p$-subgroups.
\item[\emph{(b)}] \emph{(\cite[Theorem B]{CD}).} $G$ is $\{p,q\}$-soluble of $\{p,q\}$-length at most 1, and both the Sylow $p$-subgroups and Sylow $q$-subgroups of $G$ are abelian.
\end{itemize}
\end{theorem}

It is worth mentioning that the proof of the $\{p,q\}$-solubility of $G$ uses the CFSG via the next theorem.

\begin{theorem}\emph{(\cite[Theorem 9]{CD}).}
\label{fitting}
Let $G$ be a finite group with $\fit{G}=1$. Then $\Delta(G)$ is complete.
\end{theorem}

In particular, simple groups have complete class-size prime graph. 

As it can be perceived, non-adjacency between vertices of $\Delta(G)$ highly restricts the group structure. Indeed, the most extreme case was analysed in \cite{BHM} by Bertram, Herzog and Mann, where they characterised the structure of those groups $G$ whose \emph{class-size common divisor graph} $\Gamma(G)$ is disconnected; the vertices of $\Gamma(G)$ are the non-central conjugacy classes of $G$, and two classes are adjacent if their sizes are not coprime. Certainly, $\Gamma(G)$ is connected if and only if $\Delta(G)$ is connected, so these authors particularly provided a characterisation of those groups whose class-size prime graph is disconnected. Later on, this fact was independently proved by Dolfi in \cite[Theorem 4]{dolfi1995} using different tools:

\begin{theorem}\emph{(\cite[Theorem 2]{BHM}).}
\label{disconnected}
Let $G$ be a finite group. Then $\Delta(G)$ is disconnected if and only if up to abelian direct factors, $G=A\rtimes B$ with $A$ and $B$ abelian, coprime Hall subgroups, and $G/\ze{G}$ is a Frobenius group.
\end{theorem}

In this last situation, $\Delta(G)$ turns out to be the union of two complete connected components (also called \emph{cliques}), which are the sets of prime divisors of the Frobenius kernel and a Frobenius complement of $G/\ze{G}$, respectively.

In parallel to this research, Chillag and Herzog examined in \cite{CH} a more restrictive situation, namely groups all of whose class sizes are prime powers. That result, which can be viewed as a consequence of the above one, can be restated in terms of the class-size prime graph as follows.

\begin{corollary}\emph{(\cite[Theorem 2 and Corollary 2.2]{CH}).}
\label{no-edges}
Let $G$ be a finite group. Then $\Delta(G)$ has no edges if and only if up to abelian direct factors:
\begin{itemize}
\setlength{\itemsep}{-1mm}
\item[\emph{(a)}] $G$ is a non-abelian $p$-group, for some prime $p$; or
\item[\emph{(b)}] $G=Q\rtimes P$ with $P$ a $p$-group, $Q$ a $q$-group ($p, q$ primes), both $P$ and $Q$ abelian, and $G/\ze{G}$ a Frobenius group.
\end{itemize}
\end{corollary}

Interestingly, that feature of $\Delta(G)$ in the disconnected case of having two complete induced subgraphs is not sparse. In fact, Alfandary proved in \cite{Alfandary95} that if $G$ is soluble and it does not contain a certain subgroup of an affine semi-linear group, then $\Delta(G)$ also possesses two complete induced subgraphs. Further, this situation turns out to hold in full generality for every finite group (see Corollary \ref{cor_bipartite}).

As it will be shown in the next section (see Theorem \ref{diam3}), the diameter of $\Delta(G)$ in the connected case is at most $3$. Casolo and Dolfi in 1996 characterised those groups for which the diameter of $\Delta(G)$ attains this maximal bound:

\begin{theorem}\emph{(\cite[Theorem 8]{CDdiam3}).}
Let $G$ be a group such that $\Delta(G)$ is connected. Then the diameter of $\Delta(G)$ is equal to $3$ if and only if $G=A\rtimes B$ where $A$ and $B$ are abelian, coprime Hall subgroups, and there exist $p\in \pi(B)$, $q\in\pi(A)$, $B_p\in\syl{p}{B}$, and $A_q\in\syl{q}{A}$ such that:
\begin{itemize}
\setlength{\itemsep}{-1mm}
\item[\emph{(a)}] $1\neq \ce{A}{B_p}\neq A$, $1\neq \ce{B}{A_q}\neq B$;
\item[\emph{(b)}]  $\ce{A}{y}\leqslant \ce{A}{B_p}$ for all $y\in B\smallsetminus \ze{G}$;
\item[\emph{(c)}] $\ce{B}{x}\leqslant \ce{B}{A_q}$ for all $x\in A\smallsetminus \ze{G}$.
\end{itemize}
\end{theorem}

In particular, in \cite[Proposition 7]{CDdiam3} it is also showed that if $G$ is non-soluble, then the diameter of $\Delta(G)$ is at most $2$. This bound is best possible, as Casolo and Dolfi noticed with the following example: take the direct product of an alternating group on $5$ letters and a metacyclic Frobenius group of order coprime to $30$.

The next research on this topic focused on refining the previously mentioned Itô’s result (Theorem \ref{ito} (a)). That theorem explores the structure of $G$ focusing on a prime $p$ and a single prime non-adjacent to it. In this spirit, Casolo and Dolfi considered in \cite[Lemma 5]{CDdiam3} a global, rather than local, perspective; that is, they analysed the structure of a soluble group $G$ via the non-adjacency of $p$ and a set of vertices $\pi$ distinct from $p$, by proving that the Hall $\pi$-subgroups of $G$ are abelian. Moreover, in the joint work \cite{CDPSincomplete} of these authors with Pacifici and Sanus, the solubility hypothesis was removed:

\begin{theorem}\emph{(\cite[Theorem C]{CDPSincomplete}).}
Let $G$ be a group, $p$ a vertex of $\Delta(G)$, and $\pi$ a subset of vertices that are non-adjacent to $p$ and different from $p$. Then G is $\pi$-soluble, and the following conclusions hold.
\begin{itemize}
\setlength{\itemsep}{-1mm}
\item[\emph{(a)}] $G$ has abelian Hall $\pi$-subgroups and $\pi$-length 1. 
\item[\emph{(b)}] The vertices in $\pi$ are pairwise adjacent.
\end{itemize}
\end{theorem}

In that paper, and still in the spirit of a global perspective, the authors also studied how the set of incomplete vertices (i.e. vertices that are not adjacent to at least one other vertex) of $\Delta(G)$ influences the group structure:

\begin{theorem}\emph{(\cite[Theorem A]{CDPSincomplete}).}
\label{meta}
Let $G$ be a group, and let $\pi_0$ be the set of incomplete vertices of $\Delta(G)$. Then $G$ admits a Hall $\pi_0$-subgroup $H$, and the following conclusions hold.
\begin{itemize}
\setlength{\itemsep}{-1mm}
\item[\emph{(a)}] $H$ is metabelian.
\item[\emph{(b)}] $H\cap G''\leqslant\fit{G}$.
\end{itemize}
\end{theorem}

Certainly, statement (b) is equivalent to saying that the Hall $\pi_0$-subgroups of $G''$ are nilpotent and subnormal in $G''$, and therefore $|G''/\fit{G''}|$ is not divisible by any prime in $\pi_0$. This fact and Burnside's $p^{\alpha}q^{\beta}$-theorem lead to the next result.

\begin{corollary}\emph{(\cite[Corollary B]{CDPSincomplete}).}
Let $G$ be a group. Then the following conclusions hold.
\begin{itemize}
\setlength{\itemsep}{-1mm}
\item[\emph{(a)}] All the primes within $\pi(G''/\fit{G''})$ are complete vertices of $\Delta(G)$.
\item[\emph{(b)}] If $\Delta(G)$ has at most two complete vertices, then $G$ is soluble.
\end{itemize}
\end{corollary}

Note that the class-size prime graph of an alternating group on 5 letters has three complete vertices, and clearly that group is non-soluble. 

In an earlier paper, Casolo \emph{et al.} addressed a more restrictive situation, namely $\Delta(G)$ has at most one complete vertex. Under this assumption, they proved in \cite[Theorem 2.4]{CDPSfew} not only the solubility of $G$, but also that $G'\fit{G}/\fit{G}$ is nilpotent, so $G$ is nilpotent-by-nilpotent-by-abelian. Further, if the prime $2$ is not a complete vertex of $\Delta(G)$, then $G$ is in fact nilpotent-by-metabelian, and they showed that this description is, in some sense, best possible (see \cite[Remark 2.5 and Example 2.6]{CDPSfew}).

In the most extreme case when no vertex of $\Delta(G)$ is complete, by Theorem \ref{meta} we get that $G$ is metabelian, and the following stronger conclusion was obtained by the previous authors:

\begin{theorem}\emph{(\cite[Theorem C]{CDPSfew}).}
Let $G$ be a group, and assume that no vertex of $\Delta(G)$ is complete. Then, up to abelian direct factors, $G = K \rtimes H$ where $K$ and $H$ are abelian, coprime Hall subgroups. Moreover, $K = G'$, $K \cap \ze{G} = 1$, and both $\pi(K)$ and $\pi(H)$ induce complete subgraphs in $\Delta(G)$.
\end{theorem}

Graphs that are non-complete and regular clearly satisfy the above property, and in that paper a complete characterisation of groups with such class-size prime graphs is attained. Before presenting it, we recall that the \emph{join} of two graphs $\Delta_1$ and $\Delta_2$, with disjoint vertex sets $V_1$ and $V_2$, is the graph $\Delta_1 \ast \Delta_2$ whose vertex set is $V_1\cup V_2$, and two vertices are adjacent whenever either one of them lies in $V_1$ and the other one in $V_2$, or they are vertices of the same $\Delta_i$ and they are adjacent in $\Delta_i$ for some $i\in\{1,2\}$. Besides, in that paper, the concept of a $\mathcal{D}$\emph{-group} is used for groups $G=A\rtimes B$ where $A$ and $B$ are abelian, coprime Hall subgroups, $\ze{G}\leqslant B$, and $G/\ze{G}$ is a Frobenius group.

\begin{theorem}\emph{(\cite[Theorem D]{CDPSfew}).}
\label{few}
Let $G$ be a group, and assume that $\Delta(G)$ is a non-complete regular graph of degree $d$ with $n$ vertices. Then, up to abelian direct factors, $G=G_1\times \cdots \times G_{n/2m}$, where $m=(n-1)-d$, and $G_i$ are pairwise coprime $\mathcal{D}$-groups satisfying that the orders of both the Frobenius kernel and the Frobenius complements of $G_i/\ze{G_i}$ have $m$ prime divisors.

Conversely, for a group $G$ of this kind, $\Delta(G)$ is the join of $n/2m$ copies of a graph having two complete connected components of $m$ vertices each.
\end{theorem}

As a direct consequence, we deduce that $\Delta(G)$ is a square with vertex set $V(G)=\{p,q,r,s\}$ and edge set $E(G)=\{\{p,q\},\{p,r\}, \{q,s\},\{r,s\}\}$ if and only if $G=A\times G_1\times G_2$ where $G_1$ and $G_2$ are $\mathcal{D}$-groups with $\pi(G_1)=\{p,s\}$ and $\pi(G_2)=\{q,r\}$, and $A$ is abelian.

A natural generalisation of this graph arises when one replaces each vertex by a set of vertices. With this in mind, a (simple undirected) graph is said to be a \emph{block square} if its vertex set can be partitioned into four disjoint, non-empty subsets $\pi_1,\pi_2,\pi_3,\pi_4$ such that no prime in $\pi_1$ is adjacent to any prime in $\pi_3$, no prime in $\pi_2$ is adjacent to any prime in $\pi_4$, and there exist vertices in both $\pi_1$ and $\pi_3$ that are adjacent to vertices in $\pi_2$ and in $\pi_4$.

The result below, which I proved in \cite{S}, provides a complete characterisation of those groups whose class-size prime graph is a block square.

\begin{theorem}\emph{(\cite[Theorem A]{S}).}
\label{block}
Let G be a finite group. Then $\Delta(G)$ is a block square if and only if, up to an abelian direct factor, $G=A\times B$ where $A$ and $B$ are $\mathcal{D}$-groups of coprime orders.
\end{theorem}

We close this section with recent research on how non-adjacency features of $\Delta(G)$ control the group structure. Let $\Delta$ be a graph with $n$ connected components. A vertex $v$ of $\Delta$ is said to be a \emph{cut-vertex} if the number of connected components of the graph $\Delta - v$, obtained by removing the vertex $v$ and all edges incident to $v$ from $\Delta$, is larger than $n$. 

Dolfi \emph{et al.} carried out in \cite[Theorem 3.3]{DPSScut} a complete, and somewhat technical, characterisation of the structure of a group $G$ such that $\Delta(G)$ has a cut vertex. In particular, such a group $G$ turns out to be soluble with Fitting height at most 3, and its Sylow $p$-subgroups are abelian for every prime $p$ distinct from the cut vertex (cf. \cite[Theorem A (a)]{DPSScut}). These authors also characterised in \cite[Theorem C]{DPSScut} groups whose class-size prime graph has two cut vertices (which is actually the maximal number of cut vertices that can occur, as it will be shown in the next section). As an application,
a classification of the finite groups $G$ such that $\Delta(G)$ has no cycle as an induced subgraph is achieved (cf. \cite[Corollary 3.4]{DPSScut}).

\section{Graphs that may occur as $\Delta(G)$}
\label{possible_graphs}

In this section, we put focus on which graphs $\Delta$ might satisfy that $\Delta(G)=\Delta$ for some finite group $G$. Recall that, in virtue of Theorem \ref{disconnected}, the number of connected components of $\Delta(G)$ is less than or equal to $2$ for every group $G$, having both components diameter $1$. It is natural, then, to examine whether the diameter of $\Delta(G)$ can be upper bounded in the connected case. This problem was addressed by Alfandary in \cite{Alfandary94}, where he proved the result below.

\begin{theorem}\emph{(\cite[Corollary 2.5]{Alfandary94}).}
\label{diam3}
Let $G$ be a group. If $\Delta(G)$ is connected, then its diameter is at most $3$.
\end{theorem}

Dolfi proved the same conclusion in \cite[Theorem 17]{dolfi1995}, and he provided an example which shows that this bound is best possible.

Roughly speaking, the class-size prime graph is very rich in edges. This property can be deduced from the next theorem, whose proof involves the CFSG via Theorem \ref{fitting}.

\begin{theorem}\emph{(\cite[Theorem A]{D}).}
\label{cycle}
Let $G$ be a finite group. If $p,q,r\in V(G)$ are pairwise distinct, then two among them are adjacent in $\Delta(G)$.
\end{theorem}

It is worthwhile to point out that this fact was firstly proved for soluble groups by Dolfi itself in \cite{dolfi1995}. From the above result, it follows again the bound of $3$ for the diameter of $\Delta(G)$ in the connected case, and also the feature of $\Delta(G)$ in the non-connected case of being the union of two complete subgraphs.

At this point, it is convenient to recall the definition of the complement graph $\overline{\Delta(G)}$ of the class-size prime graph: this graph has the same vertex set $V(G)$, but two vertices are adjacent in $\overline{\Delta(G)}$ if and only if they are not adjacent in $\Delta(G)$. With this in mind, Theorem \ref{cycle} can be expressed as follows: the graph $\overline{\Delta(G)}$ does not contain any cycle of length $3$, for any group $G$. Now it is natural to wonder whether there exist cycles of larger lengths within $\overline{\Delta(G)}$; and indeed any Frobenius group $G$ with abelian kernel and complement, both with orders divisible by two different primes, provides a cycle of length $4$ in $\overline{\Delta(G)}$. However, Dolfi \emph{et al.} showed in \cite{DPSSbipartite} that it is not possible to find cycles of length $5$, and in general of any odd length:

\begin{theorem}\emph{(\cite[Theorem A]{DPSSbipartite}).}
For any group $G$, the graph $\overline{\Delta(G)}$ does not contain any cycle of odd length. 
\end{theorem}

A first immediate consequence is that the class-size prime graph of a group cannot contain any pentagon as an induced subgraph. Further, the assertion of the above theorem is equivalent to the fact that $\overline{\Delta(G)}$ is always a bipartite graph, and so the following result directly follows.

\begin{corollary}\emph{(\cite[Corollaries B and C]{DPSSbipartite}).}
\label{cor_bipartite}
For every finite group $G$, the vertex set of $\Delta(G)$ can be partitioned in two subsets of pairwise adjacent vertices. In particular, if $\omega(G)$ is the maximum size of a clique in $\Delta(G)$, then it holds $|V(G)|\leq 2\omega(G)$.
\end{corollary}

Some results that have appeared in the previous section lead to characterisations of which non-complete regular graphs and which block square graphs may occur as $\Delta(G)$ for some group $G$. For instance, the following two results are consequences of Theorems \ref{few} and \ref{block}, respectively.

\begin{theorem}\emph{(\cite[Corollary E]{CDPSfew}).}
Let $\Delta$ be a non-complete regular graph. Then there exists a finite group $G$ such that $\Delta(G)=\Delta$ if and only if $\Delta$ is the join of $k$ copies of a graph that have two complete connected components of $m$ vertices each, for some positive integers $k$ and $m$.
\end{theorem}

\begin{theorem}\emph{(\cite[Corollary B]{S}).}
Let $\Delta$ be a block square graph. Then there exists a finite group $G$ such that $\Delta(G)=\Delta$ if and only if $\pi_i$ is a clique for each $1\leq i \leq 4$, and all the primes in $\pi_1\cup \pi_4$ are adjacent to all the primes in $\pi_2\cup \pi_3$.
\end{theorem}

In particular, it cannot occur a \emph{house graph} (i.e. a square graph in which two adjacent vertices are adjacent to an additional vertex) as $\Delta(G)$ for any group $G$.

Regarding graphs that possess a cut vertex, Dolfi \emph{et al.} obtained the next conclusions.

\begin{theorem}\emph{(\cite[Theorem A]{DPSScut}).}
Let $G$ be a group such that $\Delta(G)$ has a cut vertex $r$. Then the following conclusions hold:
\begin{itemize}
\setlength{\itemsep}{-1mm}
\item[\emph{(a)}] $\Delta(G) - r$ is a graph with two connected components, that are both complete graphs.
\item[\emph{(b)}] If $r$ is a complete vertex of $\Delta(G)$, then it is the unique complete vertex and the unique cut vertex of $\Delta(G)$. If $r$ is non-complete, then $\Delta(G)$ is a graph of diameter 3, and it can have at most two cut vertices.
\end{itemize}
\end{theorem}

Hence, $\Delta(G)$ can have at most two cut vertices for any group $G$, and a characterisation of the associated class-size prime graphs was also achieved in \cite{DPSScut}.

\begin{theorem}\emph{(\cite[Theorem D]{DPSScut}).}
Let $\Delta$ be a graph having a cut vertex. Then there exists a finite group $G$ such that $\Delta(G)=\Delta$ if and only if $\Delta$ is connected and its vertex set can be partitioned in two subsets of pairwise adjacent vertices.
\end{theorem}

\section{Variations on the class-size prime graph}
\label{variations}

Over the last years, several authors have investigated whether the entire data contained in $cs(G)$ is required for studying the structure of $G$. In particular, some smaller subsets of elements of the group $G$ have been considered, and hence it naturally arised the definition and the study of the prime graphs built on the associated conjugacy class sizes. It is worthwhile to remark that, in general, these graphs are not induced subgraphs of $\Delta(G)$.

\begin{center}
\subsection*{$p$-regular and $p$-singular conjugacy classes}
\end{center}

In 1972, Camina generalised Lemma \ref{vertex_out} by considering only $p$\emph{-regular} elements for a fixed prime $p$, i.e. elements of order not divisible by $p$. He proved that $p$ does not divide any conjugacy class size of each $p$-regular element of a group $G$ if and only if $G$ has a Sylow $p$-subgroup as a direct factor. Let $cs_{p'}(G)$ be the set of conjugacy class sizes of $p$-regular elements of a group $G$. Then this result can be rephrased in terms of $\Delta(cs_{p'}(G))$:

\medskip

\begin{lemma}\emph{(\cite[Lemma 1]{Camina}, rephrased).}
Let $G$ be a group, $p$ a prime, and $P\in\syl{p}{G}$. Then $p$ is not a vertex of $\Delta(cs_{p'}(G))$ if and only if $G=P\times H$ with $H\in\hall{p'}{G}$.
\end{lemma}

Lu and Zhang studied in \cite{LZ} the number of connected components of $\Delta(cs_{p'}(G))$, but assuming the $p$-solubility of the group $G$. Beltr\'an and Felipe attained best possible bounds for the diameter in \cite{BFp1} and \cite{BFp2}.

\begin{theorem}
Let $G$ be a $p$-soluble group, for a fixed prime $p$. Then:
\begin{itemize}
\setlength{\itemsep}{-1mm}
\item[\emph{(a)}] \emph{(\cite[Theorem 1]{LZ}).} The number of connected components of $\Delta(cs_{p'}(G))$ is at most $2$.
\item[\emph{(b)}] \emph{(\cite[Theorem 4]{BFp1} and \cite[Theorem A]{BFp2}).} If $\Delta(cs_{p'}(G))$ is connected, then its diameter is at most $3$; and if it is disconnected, then both connected components are complete subgraphs.
\end{itemize}
\end{theorem}

Observe that these bounds cannot be directly deduced from the fact that $\Delta(cs_{p'}(G))$ is a subgraph of $\Delta(G)$. 

In \cite{BFp3}, Beltr\'an and Felipe extended this study by analysing the structure of a $p$-soluble group $G$ with disconnected $\Delta(cs_{p'}(G))$. The authors did not achieved a complete characterisation, but they provided some necessary conditions, as well as some sufficient ones.

A characterisation of $p$-soluble groups such that all $p$-regular conjugacy classes have prime power size was addressed in \cite{BFp4} by the same authors:

\begin{theorem}\emph{(\cite[Theorem D]{BFp4}, rephrased).}
Let $G$ be a $p$-soluble group, for a fixed prime $p$. Then $\Delta(cs_{p'}(G))$ has no edges (i.e. is an empty graph) if and only if one of the following conclusions hold:
\begin{itemize}
\setlength{\itemsep}{-1mm}
\item[\emph{(a)}] $G$ has abelian $p$-complements. This occurs if and only if $cs_{p'}(G)$ only contains powers of $p$.
\item[\emph{(b)}] $G$ is nilpotent with abelian Sylow $r$-subgroups for all primes $r\notin\{p,q\}$, for some prime $q\neq p$. This happens if and only if $cs_{p'}(G)$ only contains powers of $q$.
\item[\emph{(c)}] $G=P\times H$ with $P\in\syl{p}{G}$ and $H\in\hall{p'}{G}$. Furthermore, $H=A\rtimes B$ with both $A$ and $B$ abelian, $H/\ze{H}$ is a Frobenius group, and the class sizes of $H$ are powers of either $q$ or $r$, for some primes $q$ and $r$ distinct from $p$. This happens if and only if $cs_{p'}(G)$ consists of powers of $q$ and $r$.
\end{itemize}
\end{theorem}

In particular, $\Delta(cs_{p'}(G))$ is an isolated vertex if and only if either (a) or (b) above holds. These authors extended the previous necessary condition in \cite[Theorem A]{BFpi}, by removing the $p$-solubility assumption, and considering prime power class lengths of $\pi$-elements, for a set of primes $\pi$.

Alternatively, Qian and Wang focused on the class sizes of \emph{$p$-singular elements}, which are elements whose order is divisible by $p$. Let $cs_p(G)$ denote the set of class sizes of such elements in a group $G$. Note that if $p\in \pi(\ze{G})$, then $cs(G)=cs_p(G)$. In \cite{QW} they proved the next results.

\begin{theorem}
Let $G$ be a group, $p$ a prime, and $P\in\syl{p}{G}$. 
\begin{itemize}
\setlength{\itemsep}{-1mm}
\item[\emph{(a)}] \emph{(\cite[Theorem B]{QW}, rephrased).} $p$ is not a vertex of $\Delta(cs_p(G))$ if and only if $P$ is abelian, $P\cap P^g=1$ for each $g\in G\smallsetminus \no{G}{P}$, and $\no{G}{P}/\ce{G}{P}$ acts Frobeniusly on $P$ whenever $\ce{G}{P}<\no{G}{P}$.
\item[\emph{(b)}] \emph{(\cite[Theorem C]{QW}).} If $p\notin \pi(\ze{G})$, then $\Delta(cs_p(G))$ is connected with diameter at most $3$.
\end{itemize}
\end{theorem}

\begin{center}
\subsection*{Real conjugacy classes}
\end{center}

Various authors have investigated the sizes of the so-called \emph{real conjugacy classes} of $G$, i.e. conjugacy classes $x^G$ satisfying that $\chi(x)$ is a real number for each $\chi\in \irr{G}$. This is analogous to saying that $x^g=x^{-1}$ for some $g\in G$, so $x^G=(x^{-1})^G$. Let us denote by $cs_{\mathbb{R}}(G)$ the set of sizes of the real conjugacy classes of $G$.

In 2009, Navarro, Sanus and Tiep studied the situation where all the real class sizes of $G$ are $2$-powers, which means that the set of vertices of $\Delta(cs_{\mathbb{R}}(G))$ is just the prime $2$. In virtue of \cite[Theorem C (a)]{NST}, this happens if and only if $G$ has a normal $2$-complement $K$ and all the real elements of $G$ lie in $\ce{G}{K}$. One year before, Dolfi, Navarro and Tiep characterised when $cs_{\mathbb{R}}(G)$ satisfies the ``opposite'' situation, namely all real class sizes are odd:

\medskip

\begin{theorem}\emph{(\cite[Theorem 6.1]{DNT}, rephrased).}
Let $G$ be a group and let $P\in\syl{2}{G}$. Then $2$ is not a vertex of $\Delta(cs_{\mathbb{R}}(G))$ if and only if $P$ is normal in $G$ and each real element of $P$ lies in $\ze{P}$.
\end{theorem}

For odd primes, however, the approach is quite more difficult and it depends heavily on the CFSG, which contrasts with Lemma \ref{vertex_out}. Combining results of \cite{GNT}, \cite{IN} and \cite{T}, the next theorem follows. Recall that $\pmb{\operatorname{O}}^{p'}(G)$ is the smaller normal subgroup of $G$ such that the corresponding quotient group has order not divisible by $p$.

\begin{theorem}\emph{(\cite[Lemma 2.7]{TV}, rephrased).}
Let $G$ be a group and $p$ be an odd prime. If $p=3$, assume in addition that $G$ has no composition factor isomorphic to $SL_3(2)$. If $p$ is not a vertex of $\Delta(cs_{\mathbb{R}}(G))$, then $G$ is $p$-soluble and $\pmb{\operatorname{O}}^{p'}(G)$ is soluble. Furthermore, $\pmb{\operatorname{O}}^{2'}(G)$ has a normal Sylow $p$-subgroup $P$ and $P'\leqslant \ze{\pmb{\operatorname{O}}^{2'}(G)}$.
\end{theorem}

It is also proved in \cite[Theorem 6.2]{DNT} that the number of connected components of $\Delta(cs_{\mathbb{R}}(G))$ is at most $2$. In the extreme case when $\Delta(cs_{\mathbb{R}}(G))$ is disconnected, Tong-Viet obtained the following result.

\begin{theorem}\emph{(\cite[Theorems A and B]{TV}).}
\label{real-tong-viet}
Let $G$ be a group such that $\Delta(cs_{\mathbb{R}}(G))$ is disconnected. Then $G$ is soluble, $2$ is a vertex of $\Delta(cs_{\mathbb{R}}(G))$, and one of the following conclusions holds.
\begin{itemize}
\setlength{\itemsep}{-1mm}
\item[\emph{(a)}] $G$ has a normal Sylow $2$-subgroup.
\item[\emph{(b)}] $\Delta(cs_{\mathbb{R}}(\pmb{\operatorname{O}}^{2'}(G)))$ is disconnected and $cs_{\mathbb{R}}(\pmb{\operatorname{O}}^{2'}(G))$ are either odd numbers or powers of $2$.
\end{itemize}
\end{theorem}

Actually, in the above statement (b), both connected components are complete and one of them is just the prime $2$ (see \cite[Theorem 3.5]{TV}). Concerning statement (a), Tong-Viet conjectured that it cannot occur; in other words, $\Delta(cs_{\mathbb{R}}(G))$ is connected when $G$ has a normal Sylow $2$-subgroup. In fact, he proved that conjecture whenever the real elements of the Sylow $2$-subgroup $P$ of $G$ are central in $P$ (see \cite[Lemma 3.6]{TV}).

Regarding the solubility assertion of Theorem \ref{real-tong-viet}, the same conclusion was achieved in \cite[Theorem 3.1]{DPS} by assuming a more restrictive hypothesis, namely all real class sizes of $G$ are either odd or powers of $2$. Thus, every group whose real classes all have prime power size is soluble, and its structure seems to be characterised in a current preprint due to Bonazzi (see \cite{Bo}). On the other hand, the groups whose graph $\Delta(cs_{\mathbb{R}}(G))$ consists of only one odd vertex have not yet been characterised.

\begin{center}
\subsection*{Vanishing conjugacy classes}
\end{center}

In parallel to the aforementioned research on real class sizes, some authors started to filter the class sizes of $G$ by the set $\irr{G}$ in a different manner. Concretely, an element $g\in G$ is said to be \emph{vanishing in} $G$ if $\chi(g)=0$ for some $\chi\in\irr{G}$, and so $g^G$ is called a \emph{vanishing conjugacy class}. This definition is motivated by a celebrated theorem due to Burnside, which asserts that a non-linear irreducible character always vanishes on some conjugacy class; so it is natural to analyse the columns of the character table of a group that verify that property. We refer the reader to the expository paper \cite{DPSsurvey} for a collection of results and conjectures on this topic.

Let $cs_v(G)$ be the set of vanishing class sizes of $G$. An immediate consequence of the mentioned Burnside's result is the next one: if $\Delta(cs_v(G))$ is a null graph (i.e. it has no vertices), then $G$ is abelian. Note that the converse also holds, so the commutativity of $G$ can be characterised in terms of $\Delta(cs_v(G))$.

In \cite[Theorem A]{DPSvan}, the authors studied (via the CFSG) the situation when a given prime $p$ does not divide any vanishing conjugacy class size of a group, which has the following transcription in terms of $\Delta(cs_v(G))$.

\medskip

\begin{theorem}\emph{(\cite[Theorem A]{DPSvan}, rephrased).}
\label{van_vertex}
Let $G$ be a group, and $p$ a prime. If $p$ is not a vertex of $\Delta(cs_v(G))$, then $G$ has a normal $p$-complement and abelian Sylow $p$-subgroups.
\end{theorem}

It is significant to mention that the same assertion remains true if we consider only the class sizes of elements that are zeros of some (ordinary) irreducible character of $G$ lying in the \emph{principal $p$-block}, i.e. characters $\chi\in\irr{G}$ such that the sum of $\chi(g)$, for all $p$-regular elements $g\in G$, is non-zero (cf. \cite[Theorem 6.3]{DPSsurvey}).

Certainly, the vertex set of $\Delta(cs_v(G))$ can be smaller than that of $\Delta(G)$: it is enough to consider a symmetric group on $3$ letters. Surprisingly enough, if $G$ possesses a non-abelian minimal normal subgroup, then Bianchi \emph{et al.} proved that these two vertex sets turn out to be equal:

\begin{theorem}\emph{(\cite[Proposition]{BBCP}).}
Let $G$ be a group that has a non-abelian minimal normal subgroup. Then the vertex sets of $\Delta(G)$ and $\Delta(cs_v(G))$ are equal.
\end{theorem}

In the spirit of \cite{CD}, and under the assumption that $G$ has a non-abelian minimal normal subgroup, these authors also obtained the next ``vanishing versions'' of Theorems \ref{ito} (b) and Theorem \ref{fitting}:

\begin{theorem}\emph{(\cite[Theorem A]{BBCP}).}
\label{vanishing_A}
Let $G$ be a group having a non-abelian minimal normal subgroup. If $\{p,q\}$ are vertices of $\Delta(cs_v(G))$ that are non-adjacent, then $G$ is $\{p,q\}$-soluble.
\end{theorem}

\begin{theorem}\emph{(\cite[Theorem B]{BBCP}).}
Let $G$ be a group with $\fit{G}=1$. Then $\Delta(cs_v(G))$ is a complete graph with vertex set $\pi(G)$.
\end{theorem}

An immediate $p$-solubility criterion follows from Theorems \ref{van_vertex} and \ref{vanishing_A}: if $p$ is not a complete vertex of $\Delta(cs_v(G))$, and $G$ has a non-abelian minimal normal subgroup, then $G$ is $p$-soluble (cf. \cite[Corollary]{BBCP}). Besides, as noticed in \cite{BBCP}, the assumption about the existence of a non-abelian minimal normal subgroup in Theorem \ref{vanishing_A} is fundamental. 

All the previous results might suggest that, under the assumption that $G$ has a non-abelian minimal normal subgroup, the graphs $\Delta(G)$ and $\Delta(cs_v(G))$ actually coincide, which is an open question.

In 2019, Felipe, Martínez-Pastor and myself analysed the influence of class sizes of vanishing elements of prime power order on the structure of certain factorised groups. In particular, from \cite[Corollary 4]{FMOvan} it can be obtained some information about groups satisfying that the vanishing class sizes are all prime powers, i.e. groups $G$ such that $\Delta(cs_v(G))$ is an empty graph.

\begin{corollary}
Let $G$ be a group. If $\Delta(cs_v(G))$ has no edges, then $G/\fit{G}$ is abelian.
\end{corollary}

This result was also proved independently by Robati and Hafezieh-Balaman in \cite[Corollary 2.3]{RH}. In that paper, it is additionally addressed the particular situation in which the lenghts of the vanishing conjugacy classes of a group are all powers of a fixed prime $p$; in other words, it is addressed the feature of $\Delta(cs_v(G))$ of being a single vertex $p$ (see \cite[Theorem 1.1]{RH}).

Regarding the number of connected components of $\Delta(cs_v(G))$, it has recently been demonstrated in \cite[Theorem 2.2]{RH2} that there are at most two, as it occurs in the ordinary graph $\Delta(cs(G))$. In that paper, the authors also obtained the following theorem concerning the structure of soluble groups $G$ such that $\Delta(cs_v(G))$ has exactly two connected components.

\begin{theorem}\emph{(\cite[Theorem 1.1]{RH2}).}
Let $G$ be a soluble group, and let $\pi$ be the vertex set of $\Delta(cs_v(G))$. Suppose that $\Delta(cs_v(G))$ is non-connected. Then $G=N\rtimes H$ where $N$ is the normal Hall $\pi$-subgroup of $G$ and $H$ is an abelian complement of $N$. Moreover, $N$ is a nilpotent subgroup if and only if $N$ is abelian, $G/\ze{G}$ is a Frobenius group, and there exists some element in $N$ which is vanishing in $G$.
\end{theorem}

The authors of \cite{EGA} have further filtered the set of class sizes of $G$ by considering only the \emph{SM-vanishing conjugacy classes}, i.e. conjugacy classes that are zeros of the \emph{strongly monolithic} characters of $G$. It is worth mentioning that some of the above situations also hold for the corresponding prime graph built on this last subset of class sizes of $G$ (cf. \cite{EGA}).

\begin{center}
\subsection*{Normal subgroups}
\end{center}

Beltr\'an, Felipe and Melchor considered in \cite{BFM} the set $cs_G(N)$ of conjugacy class sizes in $G$ of elements of a normal subgroup $N$ of $G$, and they introduced the associated prime graph $\Delta(cs_G(N))$. In that paper, they obtained the next result.

\bigskip

\begin{theorem}\emph{(\cite[Theorems C and D]{BFM}).}
Let $N$ be a normal subgroup of a group $G$. Then the number of connected components of $\Delta(cs_G(N))$ is at most $2$. Moreover, if $\Delta(cs_G(N))$ is connected, then its diameter is at most $3$; and if $\Delta(cs_G(N))$ is disconnected, then each connected component is a complete graph.
\end{theorem}

Additionally, in the disconnected case, the authors also obtained strong restrictions on the structure of the normal subgroup $N$:

\begin{theorem}\emph{(\cite[Theorem E]{BFM}).}
Let $N$ be a normal subgroup of a group $G$. If $\Delta(cs_G(N))$ is disconnected, then either $N=P\times A$ with $P$ a $p$-group and $A\leqslant \ze{G}$, or $N=K\rtimes H$ where both $K$ and $H$ are abelian subgroups and $N/\ze{N}$ a Frobenius group.
\end{theorem}

Note that in the previous result it appears a structure case for $N$ that differs with the general case of $N=G$ in Theorem \ref{disconnected}.

\begin{center}
\subsection*{Products of groups}
\end{center}

Another variation arises when one consider a factorised group and the conjugacy class sizes of only the elements that lie in the factors. More concretely, if $G=AB$ is the product of two subgroups $A$ and $B$, then we can define $cs_G(A\cup B)$ as the set of class sizes of the elements that belong to $A\cup B$. Felipe \emph{et al.} proved in \cite[Corollary 1]{FKMS} that, in a factorised group $G=AB$, all the prime power order elements in $A\cup B$ have class sizes not divisible by a fixed prime $p$ if and only if $G$ has a central Sylow $p$-subgroup.  Therefore, the next consequence directly follows.

\bigskip

\begin{corollary}
Let the group $G = AB$ be the product of subgroups $A$ and $B$, and let $p$ be a prime. Then $p$ is not a vertex of $\Delta(cs_G(A\cup B))$ if and only if $G$ has a central Sylow $p$-subgroup.
\end{corollary}

The proof of \cite[Corollary 1]{FKMS} uses the CFSG, in contrast to the elementary proof of Lemma \ref{vertex_out}.

In \cite{FMS}, Felipe, Martínez-Pastor and myself introduced the concept of \emph{Baer factorisations} $G=AB$, that is, products of groups such that the class sizes in $G$ of the prime power order elements in $A\cup B$ are also prime powers. This is a generalisation of the so-called \emph{Baer groups} (see \cite{B}). In \cite[Corollary 5]{FMS}, we characterised the structure of $G$ in the particular situation when those class sizes are precisely powers of a fixed prime $p$. In particular, this holds when $cs_G(A\cup B)$ only contains powers of $p$, and so the next result follows.

\begin{corollary}
Let the group $G = AB$ be the product of subgroups $A$ and $B$. Then $\Delta(cs_G(A\cup B))$ is an isolated vertex $p$ if and only if $G$ has a central $p$-complement.
\end{corollary}

More generally, we obtained some properties that are verified by Baer factorisations, and thus we can derive the same conclusions for a factorised group $G=AB$ such that $\Delta(cs_G(A\cup B))$ has no edges. For instance, we proved that $G/\fit{G}$ is abelian, and both $A$ and $B$ are Baer groups (cf. \cite{FMS}). Indeed, using \cite[Proposition D]{FMS}, it can be verified that if $\Delta(cs_G(A\cup B))$ has no edges, then both $\Delta(A)$ and $\Delta(B)$ have no edges (although in general they are not subgraphs of $\Delta(G)$), and therefore by Corollary \ref{no-edges} we can control the structure of both subgroups.

As demonstrated in the previous subsections, if some information on $cs(G)$ is not considered, then it is not always possible to get analogous results to the ordinary case. As an example, Itô's result (Theorem \ref{ito} (a)) states that if $\{p,q\}$ are non-adjacent vertices of $\Delta(G)$, then $G$ has either a normal $p$-complement or a normal $q$-complement. It is natural to wonder whether an analogous statement can be obtained when we consider this property on $\Delta(cs_G(A\cup B))$. Nevertheless, if $G=A\times B$ is the direct product of a symmetric group $A$ on $3$ letters and a dihedral group $B$ of order $10$, then $cs_G(A\cup B)=\{2,3,5\}$; so $3$ and $5$ are non-adjacent vertices of $\Delta(cs_G(A\cup B))$, and $G$ has neither a normal $3$-complement nor a normal $5$-complement. But $G$ has $\{3,5\}$-length $1$, and both the Sylow $3$-subgroup and $5$-subgroup are abelian, so it is an open question whether the claims in Theorem \ref{ito} (b) hold for factorised groups too. 

Regarding the number of connected components of $\Delta(cs_G(A\cup B))$, following the techniques in \cite{BFM} it is not difficult to prove that it can be at most $4$. But the structure of $A$ and $B$ (and possibly that of $G$) in the disconnected case, and the diameter of $\Delta(cs_G(A\cup B))$ in the connected case are problems not addressed at the moment.




\end{document}